\documentstyle[12pt]{article}
\def\Bbb R{{\rm \bf R}}
\def\proclaim#1{\vskip2mm{\bf #1}\em}
\def\endproclaim{\em \vskip2mm}
\def\tag#1{\eqno(#1)}
\def\text{\mbox}

\begin{document}

\title {Mittag-Leffler's function, Vekua transform\\
and an inverse obstacle scattering problem}
\author{Masaru IKEHATA\\
Department of Mathematics, Graduate School of Engineering\\
Gunma University, Kiryu 376-8515, JAPAN}
\date{ }
\maketitle
\begin{abstract}
This paper studies a prototype of inverse obstacle scattering problems
whose governing equation is the Helmholtz equation in two
dimensions. An explicit method to extract information about the
location and shape of unknown obstacles from the far field
operator with a fixed wave number is given.
The method is based on: an explicit construction of a modification
of Mittag-Leffler's function via the Vekua transform and the study
of the asymptotic behaviour; an explicit density in the Herglotz
wave function that approximates the modification of
Mittag-Leffler's function in the bounded domain surrounding
unknown obstacles; a system of inequalities derived from Kirsch's
factorization formula of the far field operator.
Then an indicator function which can be calculated from the far field operator acting
on the density is introduced.
It is shown that the asymptotic behaviour of the indicator function
yields information about the {\it visible part} of the exterior of
the obstacles.

\noindent
AMS: 35R30


\end{abstract}


\section{Introduction}

This paper is concerned with developing an {\it explicit}
analytical method for so-called inverse obstacle scattering
problems at a fixed wave number.  For the purpose we consider an
inverse obstacle scattering problem in two dimensions in which the
governing equation is given by the Helmholtz equation. The problem
is to extract information about the location and shape of unknown
sound-hard obstacles $D$ embedded in a medium with constant
acoustic speed and density, from the leading term of the
asymptotic expansion of the reflected wave $w$ at infinity which
is caused by an incident plane wave $e^{ikx\cdot d}$ for
infinitely many incident directions $d\in S^1$ and a fixed wave
number $k>0$. This is a prototype of several inverse obstacle
scattering problems of acoustic wave.

More precisely we assume
that: $D\subset\Bbb R^2$ is open and $\Bbb R^2\setminus\overline
D$ is connected; $\partial D$ is Lipschitz. The reflected wave $w$
is the unique solution of the scattering problem:
$$\begin{array}{c}
\displaystyle
\triangle w+k^2w=0\,\,\text{in}\,\Bbb R^2\setminus\overline D,\\
\\
\displaystyle
\frac{\partial w}{\partial\nu}=-\frac{\partial}{\partial\nu}e^{ikx\cdot d}\,\,
\text{on}\,\partial D,\\
\\
\displaystyle
\lim_{r\longrightarrow\infty}\sqrt{r}\left(\frac{\partial w}{\partial r}-ikw\right)=0
\end{array}
$$
where $\nu$ is the unit outward normal relative to $\partial D$ and $r=\vert x\vert$.
This last condition is called the Sommerfeld radiation condition.

\noindent
It is well known that given $\varphi\in S^1$ the value $w(r\varphi)$ as $r\longrightarrow\infty$ has the following
form:
$$\displaystyle
w(r\varphi)=\frac{e^{ikr}}{\sqrt{r}}F_D(\varphi;d,k)+O\left(\frac{1}{r^{3/2}}\right).
$$
The coefficient $F_D(\varphi;d,k)$ is called the far field pattern of $w$.

The operator $F: L^2(S^1)\longrightarrow L^2(S^1)$ given by the formula
$$
F_Dg(\varphi)=\int_{S^1}F_D(\varphi;d,k)g(d)d\sigma(d),\,\,g\in L^2(S^1),
$$
is called the far field operator.
It is welll known that the far field operator for a {\it fixed} $k$ uniquelly dertmines the obstacles \cite{KK}.
In this paper we consider: how to extract information about the location and shape of $D$ from the far field
operator or its partial knowledge at fixed $k$.

In \cite{Ik2, Ik3} the author established the {\it reconstruction
formula} of $D$ itself from the far field operator. The formula
consists of two parts: a relationship between a suitable
Dirichlet-to-Neumann map on the boundary of a domain that contains
$\overline D$ and the far field operator; application of the {\it
probe method} introduced by the author \cite{Ik1} to the
Dirichlet-to-Neumann map. 

In \cite{P2} Potthast gave a
reconstruction procedure that he calls the {\it singular sources
method}. The method yielded a way from the far field operator to a
scattered field outside unknown obstacles which was exerted by a
point source outside the obstacles and blows up on the boundary of
the obstacles.

Kirsch gave two types of reconstruction formulae of $D$ in \cite{K, K2}.  
The idea behind the formulae is called the {\it factorization
method} since the formulae are based on a {\it factorization formula} of the far field operator.  
In particular, in \cite{K2} he made use of the quadratic form
$$\displaystyle
(F_Dg,g)_{L^2(S^1)}
=\int_{S^1}(F_Dg)(\varphi)\overline g(\varphi) d\sigma(\varphi)
$$
acting on densities $g\in L^2(S^1)$ to introduce his indicator function.
It is defined by 
$$\displaystyle
K(x)=\inf\{\vert(F_Dg,g)\vert\,\vert\,g\in L^2(S^1),\,(\Phi_x,g)_{L^2(S^1)}=1\}
$$
where $\Phi_x(\varphi)=e^{-ikx\cdot\varphi},\,\,\varphi\in S^1$.  
He established the one line formula
$$\displaystyle
D=\{x\in\Bbb R^2\,\vert\,K(x)>0\}.
$$
For applications of his method to obstacles with other boundary conditions see \cite{G, GK}.

In \cite{Ik4} in three dimensions the author gave an extraction formula of the
convex hull of $D$ with a constraint on the Gaussian curvature of
$\partial D$ from a
Dirichlet-to-Neumann map calculated from the far field operator.
See also \cite{Ik6} for the sound-soft obstacles.
It is an application of the {\it enclosure method} introduced by the
author \cite{Ik6} and based on the asymptotic behaviour of the
function
$$\displaystyle
v=e^{x\cdot(\tau\,\vartheta+i\sqrt{\tau^2+k^2}\,\vartheta^{\perp})}
$$
having large parameter $\tau$ where both $\vartheta$ and
$\vartheta^{\perp}$ are unit vectors and perpendicular to each
other.  This function satisfies the Helmholtz equation $\triangle
v+k^2v=0$ in the whole space and divides the whole space into two
parts: if $x\cdot\vartheta>t$, then $e^{-\tau t}\vert
v\vert\longrightarrow \infty$ as $\tau\longrightarrow\infty$; if
$x\cdot\vartheta<t$, then $e^{-\tau t}\vert v\vert\longrightarrow
0$ as $\tau\longrightarrow\infty$. The indicator function
introduced in \cite{Ik4} tells us whether given $t$ the half space
$x\cdot\vartheta>t$ touches unknown obstacles.

The aim of this paper is to generalize this result by introducing
another indicator function which is given by the form $(F_Dg,g)$
acting on {\it explicit} densities $g$ on $S^1$ and tells us whether a given cone
touches unknown obstacles.

\subsection{Statement of the main result and a corollary}

In this paper we identify the point $\vartheta=(\vartheta_1,\vartheta_2)\in S^1$
with the complex number $\vartheta_1+i\vartheta_2$ and denote it by the same symbol $\vartheta$.

{\bf\noindent Definition 1.1.}
Given $n\ge 1$, $N=1,\cdots$, $s>0$ and $(y,\omega)\in B_R\times S^1$ define the indicator function
$$\displaystyle
I_{(y,\,\omega)}^{1/n}(s)_N=\left(F_Dg_{(y,\,\omega)}^{1/n}(\,\cdot\,;s,k)_N,
g_{(y,\,\omega)}^{1/n}(\,\cdot\,;s,k)_N\right)_{L^2(S^1)}
$$
where
$$\displaystyle
g_{(y,\,\omega)}^{1/n}(\varphi;s,k)_N
=\frac{\displaystyle
e^{-iky\cdot\varphi}}{2\pi}\sum_{m=0}^{n\,N}\frac{\Gamma(m+1)}
{\displaystyle\Gamma\left(\frac{m}{n}+1\right)}
\left(\frac{s\overline\omega\varphi}{ik}\right)^m,\,\,\varphi\in S^1.
\tag {1.1}
$$
Let $0<\alpha\le 1$.
Let $C_{y}(\omega,\pi\alpha/2)$ denote the interior of the cone about $\omega$ of opening angle
$\pi\alpha/2$ with vertex at $y$:
$$\displaystyle
C_{y}(\omega,\pi\alpha/2)
=\{x\in\Bbb R^2\,\vert\,
(x-y)\cdot\omega>\vert x-y\vert\cos\,(\pi\alpha/2)\}.
$$

\noindent
The following theorem is the main result of this paper.

\proclaim{\noindent Theorem 1.1.}
Let $k^2$ be not a Neumann eigenvalue of $-\triangle$ in $D$.
Assume that $\partial D$ is $C^2$ and that $\overline D$
is contained in the open disc $B_R$ radius $R$ centered at the origin.
Let $\gamma_0$ be the unique positive solution of the equation
$\displaystyle\log t+t/e=0$.
Let $\gamma$ satisfy $0<\gamma<\gamma_0$.  Let $n\ge 1$.   Let $\{s(N)\}_{N=1,\cdots}$ be an arbitrary
sequence of positive numbers satisfying, as $N\longrightarrow\infty$
$$\displaystyle
(Rs(N))^n=\frac{\gamma}{e}N+O(1).
$$
Then, given $(y,\omega)\in B_R\times S^1$ we have:

if $\overline{C_{y}(\omega,\pi/2n)}\cap\overline D=\emptyset$, then
$\lim_{N\longrightarrow\infty}\vert I_{(y,\,\omega)}^{1/n}(s(N))_N\vert=0$;

if $C_{y}(\omega,\pi/2n)\cap D\not=\emptyset$, then
$\lim_{N\longrightarrow\infty}\vert I_{(y,\,\omega)}^{1/n}(s(N))_N\vert=\infty$.

\endproclaim

Theorem 1.1 is a direct consequence of two lemmas below.

\proclaim{\noindent Lemma 1.1.}
There exists a positive constant $C$ such that, for all
$g\in L^2(S^1)$
$$\displaystyle
C^{-1}\Vert Hg\vert_{\partial D}\Vert_{H^{1/2}(\partial D)}^2\le\vert(F_Dg,g)_{L^2(S^1)}\vert
\le C\Vert Hg\vert_{\partial D}\Vert_{H^{1/2}(\partial D)}^2
$$
where
$$\displaystyle
Hg(x)=\int_{S^1}e^{ikx\cdot\varphi}g(\varphi)d\sigma(\varphi),\,\,x\in\Bbb R^2
$$
and is called the Herglotz wave function with the density $g$.

\endproclaim

\proclaim{\noindent Lemma 1.2.}
Given $(y,\omega)\in B_R\times S^1$ we have:

if $\overline{C_{y}(\omega,\pi/2n)}\cap\overline D=\emptyset$, then
$\lim_{N\longrightarrow\infty}\Vert Hg^{1/n}_{(y,\,\omega)}(\,\cdot\,;s(N),k,\omega)_N\Vert_{H^1(D)}=0$;

if $C_{y}(\omega,\pi/2n)\cap D\not=\emptyset$, then
$\lim_{N\longrightarrow\infty}\Vert Hg^{1/n}_{(y,\,\omega)}(\,\cdot\,;s(N),k,\omega)_N\vert_{\partial D}
\Vert_{L^2(\partial D)}=\infty$.

\endproclaim

Lemma 1.1 has been pointed out in \cite{G}.  It is a corollary of a factorization formula in \cite{K}
and a coerciveness of an operator in the formula.
A known proof of the coerciveness is given by a contradiction argument (cf. Lemma 4.2 in \cite{GK}) and
therefore not direct.  It seems that at the present time, there is no {\it direct proof} of
this fact.

Lemma 1.2 follows from Corollary 2.1 and (3.16) in Sections 2 and 3, respectively.

So from Theorem 1.1 what information about unknown obstacles was extracted?
To answer precisely we formulate the {\it visible part} of $B_R\setminus\overline D$.

{\bf\noindent Definition 1.2.}
We say that a point $y$ in $B_R\setminus\overline D$ is {\it visible} if the point $y$ can be connected with infinity
by a straight line that started at $y$ and goes to infinity without intersecting $\overline D$.
We denote by $V(B_R\setminus\overline D)$ the set of all points in $B_R\setminus\overline D$ that are visible.
We call this set the visible part of $B_R\setminus\overline D$.

\noindent It is easy to see that the point $y$ in
$B_R\setminus\overline D$ belongs to the visible part of
$B_R\setminus\overline D$ if and only if there exist $n$ and
$\omega\,\in\, S^1$ such that
$\overline{\{C_y(\omega,\,\pi\alpha/2)\}}\cap\overline
D=\emptyset$ with $\alpha=1/n$. The set $V(B_R\setminus\overline
D)$ is a non empty open set of $B_R\setminus\overline D$. If $D$
is convex, then we have $V(B_R\setminus\overline
D)=B_R\setminus\overline D$ by a separation theorem.
Needless to
say, in general this is not true, however, the complement of the
visible part of $B_R\setminus\overline D$ gives an estimation of
$\overline D$ from above.

The next theorem tells us that the asymptotic behaviour of the indicator function
$I^{1/n}_{(y,\,\omega)}(s(N))_N$ as $N\longrightarrow\infty$ for all $n$ and $(y,\,\omega)\in\,B_R\,\times S^1$
uniquely determines the visible part of $B_R\setminus\overline D$ except for a {\it thin} set.

\proclaim{\noindent Corollary 1.1.} Let $D_1$ and $D_2$ be two
obstacles such that: $k^2$ be not a Neumann eigenvalue of
$-\triangle$ in $D_j$; $\partial D_j$ is $C^2$ and that $\overline
D_j$ is contained in the open disc $B_R$ radius $R$ centered at
the origin. Assume that, for each fixed $n$ and
$(y,\omega)\in\,B_R\times\,S^1$ we have
$$\displaystyle
\lim_{N\longrightarrow\infty}\left((F_{D_1}-F_{D_2})\,g^{1/n}_{(y,\,\omega)}(\,\cdot\,;s(N),k),
g^{1/n}_{(y,\,\omega)}(\,\cdot\,;s(N),k)\right)_{L^2(S^1)}=0.
$$
Then $V(B_R\setminus\,\overline D_1)\setminus\partial D_2=V(B_R\setminus\,\overline D_2)\setminus\partial D_1$.

\endproclaim

\noindent This is derived from Theorem 1.1 as follows.  It
suffices to prove $V(B_R\setminus\,\overline D_1)\setminus\partial
D_2\subset V(B_R\setminus\overline D_2)$. Let
$y\in\,V(B_R\setminus\,\overline D_1)\setminus\partial D_2$.  Then there
exist $n$ and $\omega\,\in\, S^1$ such that
$\overline{\{C_y(\omega,\,\pi\alpha/2)\}}\cap\overline
D_1=\emptyset$ with $\alpha=1/n$. By Theorem 1.1 we have
$$\displaystyle
\lim_{N\longrightarrow\infty}\vert\left(F_{D_1}\,g^{1/n}_{(y,\,\omega)}(\,\cdot\,;s(N),k),
g^{1/n}_{(y,\,\omega)}(\,\cdot\,;s(N),k)\right)_{L^2(S^1)}\vert=0.
$$
Since $y$ does not belong to $\partial D_2$, it belongs to $D_2$ or $B_R\setminus\overline D_2$.
If $y\in\,D_2$, then from Theorem 1.1 one has
$$\displaystyle
\lim_{N\longrightarrow\infty}\vert\left(F_{D_2}\,g^{1/n}_{(y,\,\omega)}(\,\cdot\,;s(N),k),
g^{1/n}_{(y,\,\omega)}(\,\cdot\,;s(N),k)\right)_{L^2(S^1)}\vert=\infty.
\tag {1.2}
$$
A combination of this and assumption gives
$$\displaystyle
\lim_{N\longrightarrow\infty}
\vert\left(F_{D_1}\,g^{1/n}_{(y,\,\omega)}(\,\cdot\,;s(N),k),
g^{1/n}_{(y,\,\omega)}(\,\cdot\,;s(N),k)\right)_{L^2(S^1)}\vert=\infty.
$$
This is a contradiction.  So $y$ has to be in the set $B_R\setminus\overline D_2$.
If $\displaystyle \{C_y(\omega,\,\pi\alpha/2)\}\cap D_2\not=\emptyset$,
then, from Theorem 1.1 we again obtain (1.2) and the same contradiction as above.
Thus $\displaystyle \{C_y(\omega,\,\pi\alpha/2)\}\cap D_2=\emptyset$. Therefore if one chooses a larger
$n'$ than $n$, then one gets $\displaystyle\overline{\{C_y(\omega,\,\pi\alpha'/2)\}}\cap\overline D_2=\emptyset$
with $\alpha'=1/n'$.  This means that $y\in\,V(B_R\setminus\overline D_2)$.

\subsection{A brief explanation of the idea}

Here we give a brief explanation of the origin of the density $g_{(y,\,\omega)}^{1/n}(\,\cdot\,;s,k)_N$.

Finding the density is closely related to Mittag-Leffler's function
$E_{\alpha}(z)$ which is an entire function and defined by the formula
$$
E_{\alpha}(z)=\sum_{m=0}^{\infty}\frac{\displaystyle z^m}
{\displaystyle\Gamma(\alpha m+1)}
$$
where $\alpha$ is a parameter and satisfies $0<\alpha\le 1$.  The
function $E_{\alpha}(\tau(x_1+ix_2))$ of independent variables
$x=(x_1,x_2)$ with parameter $0<\tau<\infty$ is harmonic in
the whole plane.  This function divides the whole plane into two
parts as $\tau\longrightarrow\infty$: in a sector it is
exponentially growing; outside the sector decaying algebraically.
In \cite{Ik3} we applied this property of the harmonic function to
an inverse boundary value problem for an elliptic equation
$\nabla\cdot\gamma\nabla u=0$ with a discontinues coefficient
$\gamma$ which is a special, however, very important version of
the Calder\'on problem \cite{C} and a continuum model of
electrical impedance tomography.

In Section 2 we modify this harmonic function
by using the Vekua transform \cite{V,V2}(see also \cite{BG, CD})
which transforms given solution of
the Laplace equation in $\Bbb R^2$ into that of the Helmholtz
equation $\triangle v+k^2v=0$ in $\Bbb R^2$.  Using the
solution obtained by the transform, we define a special solution with a large parameter $s>0$
of the Helmholtz equation
which is denoted by $E_{\alpha}(x;s,k,\omega)$.
In particular, the function
$\displaystyle E_{1}(x;s,k,\omega)$
is the Vekua transform of the harmonic function
$$\displaystyle
\exp\,\left\{\frac{s}{2}\,(\omega_1-i\omega_2)(x_1+ix_2)\right\}
$$
where $\omega=(\omega_1,\omega_2)\,\in S^1$.
We show that the function $E_{\alpha}(x;s,k,\omega)$
has the asymptotic behaviour as $s\longrightarrow\infty$ similar to
that of the original Mittag-Leffler's function:

if $x\in C_{0}(\omega,\pi\alpha/2)$, then $\lim_{s\longrightarrow\infty}\vert E_{\alpha}(x;s,k,\omega)\vert=\infty$;

if $x\in\Bbb R^2\setminus\overline{C_{0}(\omega,\pi\alpha/2)}$, then
$\lim_{s\longrightarrow\infty}\vert E_{\alpha}(x;s,k,\omega)\vert=0$.

In Section 3 we establish the relationship between the density
$g_{(y,\,\omega)}^{1/n}(\,\cdot\,;s,k)_N$ and the function
$E_{\alpha}(x-y;s,k,\omega)$ for $\alpha=1/n$ and $y\in\overline B_R$:
$$\displaystyle
Hg_{(y,\,\omega)}^{1/n}(\,\cdot\,;s(N),k)_N (x)
\approx E_{1/n}(x-y;s(N),k,\omega),\,\,x\in\overline B_R
$$
as $N\longrightarrow\infty$. Thus one can say that
$g_{(y,\,\omega)}^{1/n}(\,\cdot\,;s,k)_N$ and $s(N)$ are chosen in
such a way that the corresponding Herglotz wave function
approximates a modification of Mittag-Leffler's function.

It
should be pointed out that the result in \cite{Ik11} is closely
related to the construction of the density. Therein the author
considered the case when $D$ is polygonal.  This means that $D$
has the expression $D=D_1\cup\cdots\cup D_m$ with $1\le m<\infty$
where $D_1,\cdots, D_m$ are simply connected open sets, polygons
and satisfy $\overline D_j\cap\overline D_{j'}=\emptyset$ for
$j\not=j'$.

\noindent The observation data are given by $F_D(\,\cdot\,;d,k)$ for
{\it fixed} $d$ and $k$ provided we know the disc $B_R$ that contains
$\overline D$.  Using the enclosure method \cite{Ik7}, we
established a direct extraction formula of the convex hull of $D$
from the quantity
$$\displaystyle
\int_{S^1}F_D(-\varphi;d,k)g(\varphi)d\sigma(\varphi)
$$
for some explicit densities $g$ independent of $D$. The one of key points is the
choice of the densities.  Those are chosen in such a way that
$$\displaystyle
Hg(x)\approx
e^{x\cdot(\tau\,\omega+i\sqrt{\tau^2+k^2}\,\omega^{\perp})},\,\,x\in\overline B_R
$$
where $\omega^{\perp}=(-\omega_2,\omega_1)$. However, to get more
than convex hull of unknown obstacles the function in the right
hand side is not enough. In this paper we give explicitly the
desired function by using the idea of the Vekua transform and
Mittag-Leffler's function.

Finally we point out that there are other approaches with a single incident plane wave:
the {\it point source method} \cite{P}, the {\it no response test} \cite{LP},
the {\it range test} \cite{PSK} and the notion of the {\it scattering support} \cite{KS}.

\section{Modified Mittag-Leffler's function}

In this section we introduce a modification of $E_{\alpha}(\tau(x_1+ix_2))$
with $0<\tau<\infty$ that satisfies the Helmholtz equation $\triangle u+k^2 u=0$
in $\Bbb R^2$ and study its asymptotic behaviour as $\tau\longrightarrow\infty$.

\noindent
The Bessel function of order $m=0,1,\cdots$ is given by the formula
$$\displaystyle
J_m(t)=\left(\frac{t}{2}\right)^m
\sum_{n=0}^{\infty}
\frac{\displaystyle (-1)^n}{\displaystyle (m+n)!n!}\left(\frac{t}{2}\right)^{2n}.
$$

{\bf\noindent Definition 2.1.}
Let $k\ge 0$ and $0<\alpha\le 1$.  Define
$$\displaystyle
E^k_{\alpha}(x;\tau)
=\sum_{m=0}^{\infty}\frac{\displaystyle\{\tau(x_1+ix_2)\}^m}{\displaystyle\Gamma(\alpha m+1)}
\left(\frac{\displaystyle 2}{k\vert x\vert}\right)^m m!J_m(k\vert x\vert),
\,\,0<\tau<\infty.
$$
Using the well known inequality (see Ex. 9.6, p.59 of \cite{O})
$$
\vert J_m(t)\vert\le \left(\frac{t}{2}\right)^m\frac{1}{m!},\,\, t\in\Bbb R,
\tag {2.1}
$$
one knows that $E^k_{\alpha}(x;\tau)$ is well defined and satisfies
$\displaystyle\vert E^k_{\alpha}(x;\tau)\vert\le E_{\alpha}(\tau\vert x\vert)$.

The idea behind Definition 2.1 is the following. Let
$x=(r\cos\theta,\, r\sin\theta)$.  One has
$$\displaystyle
E_{\alpha}(\tau(x_1+ix_2))=\sum_{m=0}^{\infty}\frac{\displaystyle\tau^m}
{\displaystyle\Gamma(\alpha m+1)}r^m e^{im\theta};
$$
$$\displaystyle
E_{\alpha}^k(x;\tau)=\sum_{m=0}^{\infty}\frac{\displaystyle\tau^m}
{\displaystyle\Gamma(\alpha m+1)}\left\{\left(\frac{2}{k}\right)^m\,m!J_m(kr)\right\}e^{i m\theta}.
$$
Therefore, by replacing $r^m$ in the expansion of $E_{\alpha}(\tau(x_1+ix_2))$ with
$(2/k)^m m!J_m(kr)$ one obtains $E_{\alpha}^k(x;\tau)$.
From (2.1) one knows that the absolute value
of $(2/k)^m m!J_m(kr)$ is not greater than $r^m$ and
$$\displaystyle
\left(\frac{2}{k}\right)^m m!J_m(kr)\sim r^m
$$
as $r\longrightarrow 0$.  In particular, we have $E_{\alpha}^k(0;\tau)=1$.

\noindent Using the change of variables $w=\sqrt{1-t}$ and the
formula
$$\displaystyle
k\vert x\vert\int_0^1(1-w^2)^mJ_1(k\vert x\vert w)dw
=1-\left(\frac{2}{k\vert x\vert}\right)^m m!J_m(k\vert x\vert),
$$
one can see that the function $E_{\alpha}^k(x;\tau)$ has the
integral representation given by the formula
$$
E_{\alpha}^k(x;\tau)
=E_{\alpha}(\tau(x_1+ix_2))
-\frac{\displaystyle k\vert x\vert}{\displaystyle 2}
\int_0^1 E_{\alpha}(\tau t(x_1+ix_2))J_1(k\vert x\vert\sqrt{1-t})
\frac{\displaystyle dt}{\displaystyle \sqrt{1-t}}.
\tag {2.2}
$$

The integral transform
$$
u(x)=v(x)-\frac{\displaystyle k\vert x\vert}{\displaystyle 2}
\int_0^1 v(tx)J_1(k\vert x\vert\sqrt{1-t})
\frac{\displaystyle dt}{\displaystyle \sqrt{1-t}}
$$
is called the Vekua transform of the function $v(x)$ into the
function $u(x)$. This transforms given solution of the Laplace
equation in $\Bbb R^2$ into that of the Helmholtz equation
$\triangle u+k^2u=0$ in $\Bbb R^2$. The formula (2.2) says that
$E_{\alpha}^k(x;\tau)$ is the Vekua transform of
$E_{\alpha}(\tau(x_1+ix_2))$ and therefore satisfies the Helmholtz
equation.

\noindent
In this section we show that $E_{\alpha}^k(x;\tau)$ as $\tau\longrightarrow\infty$ has
the almost same asymptotic behaviour as $E_{\alpha}(\tau(x_1+ix_2))$.

\noindent
In this paper, for convenience we introduce
$$\displaystyle
\hat{J}_m(t)=\left(\frac{2}{t}\right)^m m! J_m(t).
$$
From (2.1) we have $\vert \hat{J}_m(t)\vert\le 1$.  In what follows we make use of
this inequality frequently.

Let $f(z)$ be an arbitrary entire function of independent variable
$z=x_1+ix_2$.
Let $u(x;\tau)$ denote the Vekua transform of $f(\tau(x_1+ix_2))$.
$u(x;\tau)$ takes the form
$$
u(x;\tau)=f(\tau(x_1+ix_2))-\left(\frac{k\vert x\vert}{2}\right)^2
\int_0^1f(\tau t(x_1+ix_2))\hat{J}_1(k\vert x\vert\sqrt{1-t})dt.
\tag {2.3}
$$
In this section we write $\displaystyle
C(\pi\alpha/2)=C_0((1,0), \pi\alpha/2)$.

The following is useful for the treatment of
$E_{\alpha}^k(x;\tau)$ outside the cone $C(\pi\alpha/2)$.
See Appendix for the proof.

\proclaim{\noindent Lemma 2.1.}
One can write $u(x;\tau)$ and the partial derivatives as:
$$
\displaystyle
\left(\frac{\displaystyle 2}{\displaystyle k\vert x\vert}\right)^2
\left\{u(x;\tau)-f(\tau(x_1+ix_2))\right\}
=-\hat{J}_1(k\vert x\vert)\frac{1}{\tau}\int_0^{\tau}f(w(x_1+ix_2))dw+R(x;\tau)
\tag {2.4}
$$
where $R(x;\tau)$ satisfies
$$\displaystyle
\vert R(x;\tau)\vert\le \frac{1}{2}\left(\frac{\displaystyle k\vert x\vert}{\displaystyle 2}\right)^2
\frac{1}{\displaystyle \tau^2}\int_0^{\tau}w\vert f(w(x_1+ix_2))\vert dw;
\tag {2.5}
$$
$$\begin{array}{lr}
\displaystyle
\left(\frac{\displaystyle 2}{\displaystyle k\vert x\vert}\right)^2\left\{\frac{\partial}{\partial x_j}u(x;\tau)
-\tau i^{j-1}f'(\tau(x_1+ix_2))\right\}
=
-\frac{\displaystyle i^{j-1}\hat{J}_1(k\vert x\vert)}{x_1+ix_2}f(\tau(x_1+ix_2))\\
\\
\displaystyle
+\left\{-\frac{\displaystyle(-i)^{j-1}\hat{J}_1(k\vert x\vert)}
{\displaystyle x_1-ix_2}
+(\frac{k}{2})^2x_j\hat{J}_2(k\vert x\vert)\right\}\frac{1}{\tau}
\int_0^{\tau}f(w(x_1+ix_2))dw
+R_j(x;\tau)
\end{array}
\tag {2.6}
$$
where $R_j(x;\tau)$ satisfies
$$
\begin{array}{lr}
\displaystyle
\vert R_j(x;\tau)\vert
\le\frac{1}{2}\left(\frac{\displaystyle k\vert x\vert}{\displaystyle 2}\right)^2
\times\\
\\
\displaystyle
\left\{\frac{1}{\tau^2}\int_0^{\tau}w^2\vert f'(w(x_1+ix_2))\vert dw
+\frac{2\vert x_j\vert}{\displaystyle \vert x\vert^2}
\left(2+\frac{1}{3}\left(\frac{\displaystyle k\vert x\vert}{\displaystyle 2}\right)^2\right)
\frac{1}{\tau^2}\int_0^{\tau}w\vert f(w(x_1+ix_2))\vert dw\right\}.
\end{array}
\tag {2.7}
$$

\endproclaim

Let us consider the case when $x$ is outside cone $C(\pi\alpha/2)$.
It is known that, as $\vert x\vert\longrightarrow\infty$
and $x\in\Bbb R^2\setminus\overline{C(\pi\alpha/2)}$
Mittag-Leffler's function and the partial derivatives have
the asymptotic form (see \cite{B, E, Ik9}):
$$\displaystyle
E_{\alpha}(x_1+ix_2)=-\frac{1}{\displaystyle (x_1+ix_2)}
\frac{1}{\displaystyle \Gamma(1-\alpha)}
+O\left(\frac{1}{\displaystyle\vert x\vert^2}\right)
\tag {2.8}
$$
and
$$\displaystyle
\frac{\partial}{\partial x_j}\{E_{\alpha}(x_1+ix_2)\}=-\frac{\partial}{\partial x_j}
\{\frac{1}{\displaystyle (x_1+ix_2)}\}
\frac{1}{\displaystyle \Gamma(1-\alpha)}
+O\left(\frac{1}{\vert x\vert^3}\right).
\tag {2.9}
$$
These asymptotics are valid uniformly in the region $\{x\in\Bbb
R^2\setminus\overline{C(\pi\alpha/2+\epsilon)}\, \vert\, R_0<\vert
x\vert\}$ for given $\pi-\pi\alpha/2>\epsilon>0$ and some
$R_0>>1$.

\proclaim{\noindent Proposition 2.1.}
Let $x\in\Bbb R^2\setminus\overline{C(\pi\alpha/2)}$.  We have, as $\tau\longrightarrow\infty$
$$
E_{\alpha}^k(x;\tau)
=\left(\frac{k\vert x\vert}{2}\right)^2\frac{\displaystyle\hat{J}_1(k\vert x\vert)}{\displaystyle x_1+ix_2}
\frac{1}{\displaystyle \Gamma(1-\alpha)}\frac{\log\tau}{\tau}
+O\left(\frac{1}{\tau}\right)
\tag {2.10}
$$
and
$$
\frac{\partial}{\partial x_j}E_{\alpha}^k(x;\tau)
=\frac{\partial}{\partial x_j}
\left\{\left(\frac{k\vert x\vert}{2}\right)^2\frac{\displaystyle\hat{J}_1(k\vert x\vert)}
{\displaystyle x_1+ix_2}\right\}
\frac{1}{\displaystyle \Gamma(1-\alpha)}\frac{\log\tau}{\tau}
+O\left(\frac{1}{\tau}\right).
\tag {2.11}
$$
These asymptotics are valid uniformly in the region $\{x\in\Bbb
R^2\setminus\overline{C(\pi\alpha/2+\epsilon)}\, \vert\,
R^{-1}<\vert x\vert<R\}$ for given $\pi-\pi\alpha/2>\epsilon>0$
and $R>0$.
\endproclaim

{\it\noindent Proof.}
From (2.8) and (2.9) we obtain
$$\displaystyle
\frac{1}{\tau}\int_0^{\tau}E_{\alpha}(w(x_1+ix_2))dw=
-\frac{1}{x_1+ix_2}\frac{1}{\Gamma(1-\alpha)}\frac{\log\tau}{\tau}
+O\left(\frac{1}{\tau}\right);
\tag {2.12}
$$
$$\displaystyle
\frac{1}{\tau^2}\int_0^{\tau}w\vert E_{\alpha}(w(x_1+ix_2))\vert dw=O\left(\frac{1}{\tau}\right);
\tag {2.13}
$$
$$\displaystyle
\frac{1}{\tau^2}\int_0^{\tau}w^2\vert E_{\alpha}'(w(x_1+ix_2))\vert dw=O\left(\frac{1}{\tau}\right).
\tag {2.14}
$$
Then applying Lemma 2.1 to the case when $f(z)=E_{\alpha}(z)$,
from (2.4), (2.5), (2.12) and (2.13),
we obtain (2.10).
Next from (2.6), (2.7), (2.8), (2.9), (2.13) and (2.14), we have
$$
\begin{array}{lr}
\displaystyle
\frac{\partial}{\partial x_j}E_{\alpha}^k(x;\tau)
=\tau i^{j-1}E'_{\alpha}(\tau(x_1+ix_2))
-\left(\frac{k\vert x\vert}{2}\right)^2\frac{\displaystyle i^{j-1}\hat{J}_1(k\vert x\vert)}
{x_1+ix_2}E_{\alpha}(\tau(x_1+ix_2))\\
\\
\displaystyle
+\left(\frac{k\vert x\vert}{2}\right)^2
\left\{-\frac{\displaystyle (-1)^{j-1}\hat{J}_1(k\vert x\vert)}{x_1-ix_2}
+\left(\frac{k}{2}\right)^2x_j\hat{J}_2(k\vert x\vert)\right\}
\frac{1}{\tau}\int_0^{\tau}E_{\alpha}(w(x_1+ix_2))dw
+O\left(\frac{1}{\tau}\right)\\
\\
\displaystyle
=\left(\frac{k\vert x\vert}{2}\right)^2
\left\{-\frac{\displaystyle (-1)^{j-1}\hat{J}_1(k\vert x\vert)}{x_1-ix_2}
+\left(\frac{k}{2}\right)^2x_j\hat{J}_2(k\vert x\vert)\right\}
\frac{1}{\tau}\int_0^{\tau}E_{\alpha}(w(x_1+ix_2))dw
+O\left(\frac{1}{\tau}\right).
\end{array}
$$
Then from the equation
$$
\frac{\partial}{\partial x_j}
\left\{\left(\frac{k\vert x\vert}{2}\right)^2\frac{\displaystyle\hat{J}_1(k\vert x\vert)}
{\displaystyle x_1+ix_2}\right\}
=\left(\frac{\displaystyle k\vert x\vert}{2}\right)^2
\left\{\frac{\displaystyle (-i)^{j-1}}{x_1-ix_2}\hat{J}_1(k\vert x\vert)
-\left(\frac{k}{2}\right)^2x_j\hat{J}_2(k\vert x\vert)\right\}\frac{1}{\displaystyle x_1+ix_2}
$$
and (2.12) one obtains (2.11).

$\Box$

Next consider the case when $x$ is inside the cone
$C_y(\pi\alpha/2)$.
From (2.2) we have the expression
$$\displaystyle
E_{\alpha}^k(x;\tau)
=E_{\alpha}(\tau(x_1+ix_2))
-\left(\frac{k\vert x\vert}{2}\right)^2\int_0^1E_{\alpha}(\tau t(x_1+ix_2))\hat{J}_1(k\vert x\vert\sqrt{1-t})dt
\tag {2.15}
$$
and a direct computation yields
$$\begin{array}{lr}
\displaystyle
\frac{\partial}{\partial x_j}E_{\alpha}^k(x;\tau)
\\
\\
\displaystyle
=\tau i^{j-1}
\left\{E_{\alpha}'(\tau(x_1+ix_2))
-\left(\frac{k\vert x\vert}{2}\right)^2
\int_0^1 E_{\alpha}'(\tau t(x_1+ix_2))\hat{J}_1(k\vert x\vert\sqrt{1-t})dt\right\}\\
\\
\displaystyle
+\tau i^{j-1}\left(\frac{k\vert x\vert}{2}\right)^2\int_0^1
E_{\alpha}'(\tau t(x_1+ix_2))(1-t)\hat{J}_1(k\vert x\vert\sqrt{1-t})dt\\
\\
\displaystyle
-2\left(\frac{k}{2}\right)^2x_j\int_0^1E_{\alpha}(\tau t(x_1+ix_2))\hat{J}_1(k\vert x\vert\sqrt{1-t})dt\\
\\
\displaystyle
+\left(\frac{k}{2}\right)^4\vert x\vert^2 x_j
\int_0^1E_{\alpha}(\tau t(x_1+ix_2))(1-t)\hat{J}_2(k\vert x\vert\sqrt{1-t})dt.
\end{array}
\tag {2.16}
$$

It is known that there exists a positive constant $C$ such that,
for all $z\in\overline{C(\pi\alpha/2)}\setminus\{0\}$
the estimates
$$\displaystyle
\vert E_{\alpha}(z)-
\frac{1}{\alpha}e^{\displaystyle z^{1/\alpha}}\vert\le\frac{C}{\displaystyle (1+\vert z\vert^2)^{1/2}}
\tag {2.17}
$$
and
$$\displaystyle
\vert\frac{d}{dz}
\left\{E_{\alpha}(z)-\frac{1}{\alpha}
e^{\displaystyle z^{1/\alpha}}\right\}\vert\le\frac{C}{\displaystyle 1+\vert z\vert^2},
\tag {2.18}
$$
are valid (see \cite{E, Ik9}).  Note that in (2.18) there is no restriction on $z$ in a neighbourhood of $0$.
This is because of $0<\alpha\le 1$.

\proclaim{\noindent Proposition 2.2.}
Given $R>0$ and $\epsilon>0$ let $x$ satisfy $R^{-1}\le\vert x\vert\le R$ and
$\text{Re}\,(x_1+ix_2)^{1/\alpha}\ge \epsilon$.
Then, as $\tau\longrightarrow\infty$ we have two formulae:
$$\displaystyle
E_{\alpha}^{k}(x;\tau)
=\frac{1}{\alpha}e^{\displaystyle \tau^{1/\alpha}(x_1+ix_2)^{1/\alpha}}
\left\{1+O\left(\frac{1}{\displaystyle \tau^{1/\alpha}}\right)\right\};
\tag {2.19}
$$
$$\displaystyle
\frac{\partial}{\partial x_j}E_{\alpha}^k(x;\tau)
=\frac{\partial}{\partial x_j}
\left\{\frac{1}{\alpha}e^{\displaystyle \tau^{1/\alpha}(x_1+ix_2)^{1/\alpha}}\right\}
\left\{1+O\left(\frac{1}{\displaystyle \tau^{1/\alpha}}\right)\right\}.
\tag {2.20}
$$
\endproclaim
{\it\noindent Proof.}
Since (2.19) is a direct consequence of (2.15), (2.17) and the estimate
$$\displaystyle
\int_0^1 E_{\alpha}(\tau t(x_1+ix_2))\hat{J}_1(k\vert x\vert\sqrt{1-t})dt
=\frac{1}{\alpha}e^{\displaystyle \tau^{1/\alpha}(x_1+ix_2)^{1/\alpha}}
O\left(\frac{1}{\displaystyle \tau^{1/\alpha}}\right),
\tag {2.21}
$$
first we give the proof of (2.21).

\noindent
Write
$$\begin{array}{c}
\displaystyle
\int_0^1 E_{\alpha}(\tau t(x_1+ix_2))\hat{J}_1(k\vert x\vert\sqrt{1-t})dt
=\int_0^1\frac{1}{\alpha}e^{\displaystyle\tau^{1/\alpha}t^{1/\alpha}(x_1+ix_2)^{1/\alpha}}
\hat{J}_1(k\vert x\vert\sqrt{1-t})dt\\
\\
\displaystyle
+\int_0^1\left\{E_{\alpha}(\tau t(x_1+ix_2))-\frac{1}{\alpha}e^{\displaystyle\tau^{1/\alpha}
t^{1/\alpha}(x_1+ix_2)^{1/\alpha}}\right\}
\hat{J}_1(k\vert x\vert\sqrt{1-t})dt\\
\\
\displaystyle
\equiv I+II.
\end{array}
$$
Then, from (2.17) we have
$$
\displaystyle
\vert II\vert
\le C\int_0^1\frac{dt}{\displaystyle (1+\tau^2\vert x\vert^2)^{1/2}}
=\frac{C}{\tau\vert x\vert}\int_0^{\tau\vert x\vert}\frac{ds}{\displaystyle (1+s^2)^{1/2}}
=O\left(\frac{\displaystyle \log\,\tau}{\tau}\right)
\tag {2.22}
$$
provided $R^{-1}\le\vert x\vert\le R$.

\noindent
Let $0<\delta<1$.  Write
$$\begin{array}{lr}
\displaystyle
I=\int_0^{\delta}\frac{1}{\alpha}e^{\displaystyle\tau^{1/\alpha}t^{1/\alpha}
(x_1+ix_2)^{1/\alpha}}\hat{J}_1(k\vert x\vert\sqrt{1-t})dt\\
\\
\displaystyle
+\int_{\delta}^1\frac{1}{\alpha}e^{\displaystyle\tau^{1/\alpha}t^{1/\alpha}
(x_1+ix_2)^{1/\alpha}}\hat{J}_1(k\vert x\vert\sqrt{1-t})dt\\
\\
\displaystyle
\equiv I_1+I_2.
\end{array}
\tag {2.23}
$$
Then, one has
$$
\vert I_1\vert\le\frac{\delta}{\alpha}e^{\displaystyle \tau^{1/\alpha}\delta^{1/\alpha}\text{Re}\,(x_1+ix_2)^{1/\alpha}}.
\tag {2.24}
$$
Write
$$
I_2=\frac{1}{\alpha}e^{\displaystyle\tau^{1/\alpha}(x_1+ix_2)^{1/\alpha}}
\int_{\delta}^1e^{\displaystyle -\tau^{1/\alpha}(1-t^{1/\alpha})(x_1+ix_2)^{1/\alpha}}
\hat{J}_1(k\vert x\vert\sqrt{1-t})dt.
$$
Change of a variable $1-t^{1/\alpha}=s$ yields
$$\begin{array}{c}
\displaystyle
\int_{\delta}^1e^{\displaystyle -\tau^{1/\alpha}(1-t^{1/\alpha})(x_1+ix_2)^{1/\alpha}}
\hat{J}_1(k\vert x\vert\sqrt{1-t})dt\\
\\
\displaystyle
=\alpha\int_0^{1-\delta^{1/\alpha}}e^{\displaystyle -\tau^{1/\alpha}s(x_1+ix_2)^{1/\alpha}}
\frac{\displaystyle\hat{J}_1(k\vert x\vert \sqrt{1-(1-s)^{\alpha}})}{\displaystyle (1-s)^{1-\alpha}}ds
\end{array}
$$
and this gives
$$
\begin{array}{c}
\displaystyle
\vert\int_{\delta}^1e^{\displaystyle -\tau^{1/\alpha}(1-t^{1/\alpha})(x_1+ix_2)^{1/\alpha}}
\hat{J}_1(k\vert x\vert\sqrt{1-t})dt\vert\\
\\
\displaystyle
\le\alpha\int_0^{1-\delta^{1/\alpha}}
e^{\displaystyle -\tau^{1/\alpha}s\text{Re}\,(x_1+ix_2)^{1/\alpha}}
\frac{ds}{\displaystyle (1-s)^{1-\alpha}}\\
\\
\displaystyle
\le\frac{\alpha}{\displaystyle \{1-(1-\delta^{1/\alpha})\}^{1-\alpha}}
\int_0^{1-\delta^{1/\alpha}}
e^{\displaystyle -\tau^{1/\alpha}s\text{Re}\,(x_1+ix_2)^{1/\alpha}}ds\\
\\
\displaystyle
\le
\frac{\alpha}{\displaystyle \delta^{(1-\alpha)/\alpha}\tau^{1/\alpha}\text{Re}\,(x_1+ix_2)^{1/\alpha}}.
\end{array}
$$
Therefore we obtain
$$
\vert I_2\vert
\le\frac{\displaystyle e^{\displaystyle\tau^{1/\alpha}\text{Re}\,(x_1+ix_2)^{1/\alpha}}}
{\displaystyle \delta^{(1-\alpha)/\alpha}\tau^{1/\alpha}\text{Re}\,(x_1+ix_2)^{1/\alpha}}.
\tag {2.25}
$$
From (2.23), (2.24) and (2.25) one concludes that
$$
I\alpha e^{\displaystyle -\tau^{1/\alpha}(x_1+ix_2)^{1/\alpha}}
=O\left(\frac{1}{\displaystyle \tau^{1/\alpha}}\right)
\tag {2.26}
$$
provided $\text{Re}\,(x_1+ix_2)^{1/\alpha}\ge\epsilon$.  A combination of (2.22) and (2.26) yields (2.21).

Next we prove: for $m=0, 1$
$$
\begin{array}{c}
\displaystyle
\int_0^1 E'_{\alpha}(\tau t(x_1+ix_2))(1-t)^m\hat{J}_1(k\vert x\vert\sqrt{1-t})dt\\
\\
\displaystyle
=\frac{1}{\displaystyle\alpha^2}
e^{\displaystyle \tau^{1/\alpha}(x_1+ix_2)^{1/\alpha}}\tau^{(1-\alpha)/\alpha}
(x_1+ix_2)^{(1-\alpha)/\alpha}O\left(\frac{1}{\displaystyle \tau^{1/\alpha}}\right).
\end{array}
\tag {2.27}
$$

We have
$$\begin{array}{c}
\displaystyle
\int_0^1 e^{\displaystyle \tau^{1/\alpha}t^{1/\alpha}(x_1+ix_2)^{1/\alpha}}\tau^{(1-\alpha)/\alpha}
t^{(1-\alpha)/\alpha}(x_1+ix_2)^{(1-\alpha)/\alpha}
(1-t)^m\hat{J}_1(k\vert x\vert\sqrt{1-t})dt\\
\\
\displaystyle
=e^{\displaystyle \tau^{1/\alpha}(x_1+ix_2)^{1/\alpha}}
\tau^{(1-\alpha)/\alpha}(x_1+ix_2)^{(1-\alpha)/\alpha}\\
\\
\displaystyle
\times
\int_0^1 e^{\displaystyle -\tau^{1/\alpha}(1-t^{1/\alpha})(x_1+ix_2)^{1/\alpha}}
t^{(1-\alpha)/\alpha}(1-t)^m\hat{J}_1(k\vert x\vert\sqrt{1-t})dt.
\end{array}
$$
Since
$$\begin{array}{c}
\displaystyle
\vert
\int_0^1 e^{\displaystyle -\tau^{1/\alpha}(1-t^{1/\alpha})(x_1+ix_2)^{1/\alpha}}
t^{(1-\alpha)/\alpha}(1-t)^m\hat{J}_1(k\vert x\vert\sqrt{1-t})dt\vert\\
\\
\displaystyle
\le\int_0^1 e^{\displaystyle -\tau^{1/\alpha}(1-t^{1/\alpha})\text{Re}\,(x_1+ix_2)^{1/\alpha}}
t^{(1-\alpha)/\alpha}dt=\alpha\int_0^1 e^{\displaystyle -\tau^{1/\alpha} s\text{Re}\,(x_1+ix_2)^{1/\alpha}}ds\\
\\
\displaystyle
\le\frac{\alpha}{\displaystyle \tau^{1/\alpha}\text{Re}\,(x_1+ix_2)^{1/\alpha}},
\end{array}
$$
we obtain
$$\begin{array}{c}
\displaystyle
\int_0^1 e^{\displaystyle \tau^{1/\alpha}t^{1/\alpha}(x_1+ix_2)^{1/\alpha}}\tau^{(1-\alpha)/\alpha}
t^{(1-\alpha)/\alpha}(x_1+ix_2)^{(1-\alpha)/\alpha}
(1-t)^m\hat{J}_1(k\vert x\vert\sqrt{1-t})dt\\
\\
\displaystyle
=e^{\displaystyle \tau^{1/\alpha}(x_1+ix_2)^{1/\alpha}}
\tau^{(1-\alpha)/\alpha}(x_1+ix_2)^{(1-\alpha)/\alpha}
O\left(\frac{1}{\displaystyle \tau^{1/\alpha}}\right).
\end{array}
\tag {2.28}
$$

Now from (2.18) and (2.28) one obtains (2.27).
Similarly, for $j=1, 2$ we have
$$\begin{array}{c}
\displaystyle
\int_0^1 E_{\alpha}(\tau t(x_1+ix_2))(1-t)^m\hat{J}_{m+1}(k\vert x\vert\sqrt{1-t})dt\\
\\
\displaystyle
=\tau i^{j-1}\frac{1}{\displaystyle\alpha^2}
e^{\displaystyle \tau^{1/\alpha}(x_1+ix_2)^{1/\alpha}}
\tau^{(1-\alpha)/\alpha}(x_1+ix_2)^{(1-\alpha)/\alpha}
O\left(\frac{1}{\displaystyle \tau^{1/\alpha}}\right)
\end{array}
\tag {2.29}
$$
provided $\text{Re}\,(x_1+ix_2)^{1/\alpha}\ge \epsilon$ and $R^{-1}\le\vert x\vert\le R$.
Note that this is a `rough' estimate.
Now from (2.16), (2.18), (2.28) and (2.29) we obtain (2.20).

\noindent
$\Box$

As a corollary of Propositions 2.1 and 2.2 we have immediately

\proclaim{\noindent Corollary 2.1.}
We have:

for any regular $C^2$ curve $c$ with $\overline c\subset C(\pi\alpha/2)$
$$\displaystyle
\lim_{\tau\longrightarrow\infty}
\Vert E_{\alpha}^k(\,\cdot\,;\tau)\vert_{c}\Vert_{L^2(c)};
\longrightarrow\infty
$$

for any non empty bounded open set $U$ of $\Bbb R^2$ with $\overline U\subset\Bbb R^2
\setminus\overline{C(\pi\alpha/2)}$
$$\displaystyle
\lim_{\tau\longrightarrow\infty}\Vert E_{\alpha}^k(\,\cdot\,;\tau)\Vert_{H^1(U)}
\longrightarrow 0.
$$

\endproclaim

\section{Construction of the density}

{\bf\noindent Definition 3.1.}
Given $\omega=(\omega_1,\,\omega_2)\,\in S^1$ set
$\omega^{\perp}=(-\omega_2,\,\omega_1)$.
Define the function $E_{1/n}(x;s,k,\omega)$ by the formula
$$\displaystyle
E_{1/n}(x;s,k,\omega)
=E^k_{1/n}\left((x\cdot\omega,x\cdot\omega^{\perp});\frac{s}{2}\right),\,\,s>0.
$$

From (2.2) we have already known that the function $E_{1/n}(x;s,k,\omega)$
of $x\in\Bbb R^2$ satisfies the Helmholtz equation $\triangle v+k^2v=0$ in $\Bbb R^2$.
Since $\displaystyle x\cdot(\omega+i\omega^{\perp})=(\omega_1-i\omega_2)(x_1+ix_2)$, the function
$E_{1/n}(x;s,k,\omega)$ coincides with the Vekua transform of the harmonic function
$\displaystyle E_{1/n}(s(\omega_1-i\omega_2)(x_1+ix_2)/2))$ in $\Bbb R^2$.

The aim of this section is to construct a density $g\in L^2(S^1)$
explicitly such that
$$\displaystyle
Hg(x)\approx E_{1/n}(x;s,k,\omega), \,\,x\in\overline B_{2R}.
$$

\noindent
The starting point is the following fact.

\proclaim{\noindent Proposition 3.1(\cite{Ik11}).} The Vekua
transform of the harmonic function
$$\displaystyle
e^{\displaystyle ik\overline\varphi(x_1+ix_2)/2}+e^{\displaystyle ik\varphi(x_1-ix_2)/2}-1
$$
coincides with $e^{ikx\cdot\varphi}$.

\endproclaim

Let $\Gamma$ be a non empty open subset of $S^1$.
Given $g\in L^2(S^1)$ the function
$$\displaystyle
\int_{\Gamma}\left\{e^{\displaystyle ik\overline\varphi(x_1+ix_2)/2}
+e^{\displaystyle ik\varphi(x_1-ix_2)/2}-1\right\}g(\varphi)d\sigma(\varphi)
$$
is harmonic in the whole plane.  As a corollary of Proposition 3.1 one knows that
the Vekua transform of this harmonic function coincides with the Herglotz wave function
$Hg$.

\noindent Taking account the fact mentioned above and the
definition of $E_{1/n}(x;s,k,\omega)$, it suffices to construct
$g$ in such a way that
$$
\displaystyle
\int_{\Gamma}\left\{e^{\displaystyle ik\overline\varphi(x_1+ix_2)/2}+
e^{\displaystyle ik\varphi(x_1-ix_2)/2}-1\right\}g(\varphi)d\sigma(\varphi)
\approx
E_{1/n}\left(\frac{s}{2}\overline\omega(x_1+ix_2)\right)
\tag {3.1}
$$
where $\overline\omega=\omega_1-i\omega_2$.

\noindent
Using the power series expansion of Mittag-Leffler's function, one knows that if $g$ satisfies the system of equations
$$
\displaystyle
\frac{1}{\Gamma(m+1)}\left(\frac{ik}{2}\right)^m\int_{\Gamma}\varphi^mg(\varphi)d\sigma(\varphi)
=0, m=1,\cdots
\tag {3.2}
$$
and
$$
\displaystyle
\frac{1}{\Gamma(m+1)}\left(\frac{ik}{2}\right)^m\int_{\Gamma}\overline{\varphi}^m
g(\varphi)d\sigma(\varphi)
=\frac{1}{\displaystyle\Gamma\left(\frac{m}{n}+1\right)}
\left(\displaystyle\frac{s}{2}\overline\omega\right)^m, m=0,1,\cdots,
\tag {3.3}
$$
then $g$ satisfies (3.1) exactly.
Now consider the case when $\Gamma=S^1$.
We construct $g$ in the form
$$
g(\varphi)=\sum_{m=0}^{\infty}\beta_m\varphi^m+\sum_{m=1}^{\infty}\beta_{-m}\overline\varphi^{m}.
$$
Since
$$
\displaystyle
\frac{1}{2\pi}\int_{S^1}\overline\varphi^mg(\varphi)d\sigma(\varphi)=\beta_m,\,\,
\frac{1}{2\pi}\int_{S^1}\varphi^mg(\varphi)d\sigma(\varphi)
=\beta_{-m},
$$
from (3.2) and (3.3) we get
$\displaystyle\beta_{-m}=0,\,\,m=1,2,\cdots$
and
$$
\displaystyle
\beta_{m}
=\frac{1}{2\pi}\frac{\Gamma(m+1)}{\displaystyle\Gamma\left(\frac{m}{n}+1\right)}
\left(\frac{s\overline\omega}{ik}\right)^m,\,\,m=0,1,\cdots.
$$
Then $g$ becomes
$$
\displaystyle
g(\varphi)
=\sum_{m=0}^{\infty}\frac{1}{2\pi}
\frac{\Gamma(m+1)}{\displaystyle\Gamma\left(\frac{m}{n}+1\right)}
\left(\frac{s\overline\omega\varphi}{ik}\right)^m
\tag {3.4}
$$
This is always divergent.  So we consider a
truncation of (3.4):
$$
\displaystyle
g_{N}(\varphi;s,k,\omega)
=
\frac{1}{2\pi}\sum_{m=0}^{nN}\frac{\Gamma(m+1)}
{\displaystyle
\Gamma\left(\frac{m}{n}+1\right)}\left(\frac{s\overline\omega\varphi}{ik}\right)^m
\tag {3.5}
$$
where $N=1,\cdots$.
Then one obtains
$$\begin{array}{c}
\displaystyle
\int_{S^1}\left\{e^{\displaystyle ik\overline\varphi(x_1+ix_2)/2}+e^{\displaystyle ik\varphi(x_1-ix_2)/2}-1\right\}
g_N(\varphi;s,k,\omega)d\sigma(\varphi)\\
\\
\displaystyle
=
E_{1/n}\left(\frac{s}{2}\overline\omega(x_1+ix_2)\right)
-\sum_{m>nN}\frac{1}{\displaystyle\Gamma\left(\frac{m}{n}+1\right)}
\left(\frac{s}{2}\overline\omega\right)^m(x_1+ix_2)^m.
\end{array}
\tag {3.6}
$$
This shows $g_{N}(\,\cdot\,;s,k,\omega)$ satisfies (3.1) in
this sense.  Taking the Vekua transform of the both sides of
(3.6) we obtain the equation
$$
\displaystyle
Hg_N(\,\cdot\,;s,k,\omega)(x)-E_{1/n}(x;s,k,\omega)
=
-\sum_{m>nN}\frac{\Gamma(m+1)}{\displaystyle\Gamma\left(\frac{m}{n}+1\right)}
\left(\frac{s}{k}\overline\omega\right)^mJ_m(kr)e^{im\theta}
\tag {3.7}
$$
where $x=(r\cos\,\theta,r\sin\,\theta)$.  Note that this can be checked also directly
and the equation $Hg_N(\,\cdot\,;s,k,\omega)(0)=E_{1/n}(0;s,k,\omega)=1$ holds.

For our purpose we have to consider how to choose $s$ depending on $N$.
One answer to this question is the following and it is the main result of this section.

\proclaim{\noindent Theorem 3.1.}
Let $\gamma_0$ be the unique positive solution of the equation
$\displaystyle\log t+ t/e=0$.
Let $\gamma$ satisfy $0<\gamma<\gamma_0$.  Let $\{s(N)\}_{N=1,\cdots}$ be an arbitrary
sequence of positive numbers satisfying, as $N\longrightarrow\infty$
$$\displaystyle
(Rs(N))^n=\frac{\gamma}{e}N+O(1).
$$
Then we have, as $N\longrightarrow\infty$
$$\begin{array}{c}
\displaystyle \sup_{\vert x\vert\le
2R}\vert Hg_N(\,\cdot\,;s(N),k,\omega)(x)
-E_{1/n}(x;s(N),k,\omega)\vert
\\
\\
\displaystyle
+\sup_{\vert x\vert\le 2R}\vert\nabla
\{Hg_N(\,\cdot\,;s(N),k,\omega)(x)
-E_{1/n}(x;s(N),k,\omega)\}\vert
\\
\\
\displaystyle
=O\left(N^{3/2}e^{\displaystyle N(\frac{\gamma}{e}+\log\gamma)}\right)=O(N^{-\infty}).
\end{array}
\tag {3.8}
$$
\endproclaim

{\it\noindent Proof.}
Set
$$
\displaystyle
S_N(x;s)=\sum_{m>nN}\frac{\Gamma(m+1)}{\displaystyle
\Gamma\left(\frac{m}{n}+1\right)}
\left(\frac{s}{k}\overline\omega\right)^m J_m(kr)e^{im\theta},\,\,
E^{nN}_{1/n}(z)=\sum_{m=0}^{nN}\frac{z^m}{\displaystyle\Gamma\left(\frac{m}{n}+1\right)}.
$$
Then from (2.1) we have, for all $x$ with $\vert x\vert\le 2R$
$$
\displaystyle
\vert S_N(x;s)\vert
\le \sum_{m>nN}\frac{1}{\displaystyle\Gamma\left(\frac{m}{n}+1\right)}(Rs)^m
=\{E_{1/n}(z)-E^{nN}_{1/n}(z)\}\vert_{z=Rs}.
\tag {3.9}
$$
Moreover, using the recurrence relation
$$
\displaystyle
J_{m+1}(kr)=\frac{m}{kr}J_m(kr)-J_{m}'(kr),\,\,
J_{m-1}(kr)=\frac{m}{kr}J_m(kr)+J_{m}'(kr)
$$
one has the formulae
$$\begin{array}{c}
\displaystyle
e^{i\theta}\left(\frac{\partial}{\partial r}+i\frac{1}{r}\frac{\partial}{\partial\theta}\right)
J_m(kr)e^{im\theta}
=-kJ_{m+1}(kr)e^{i(m+1)\theta};
\\
\\
\displaystyle
e^{-i\theta}\left(\frac{\partial}{\partial r}-i\frac{1}{r}\frac{\partial}{\partial\theta}\right)
J_m(kr)e^{im\theta}
=kJ_{m-1}(kr)e^{i(m-1)\theta}.
\end{array}
$$
Then the formulae
$$\begin{array}{c}
\displaystyle
\frac{\partial}{\partial x_1}
=\frac{e^{i\theta}}{2}\left(\frac{\partial}{\partial r}+i\frac{1}{r}\frac{\partial}{\partial\theta}\right)
+\frac{e^{-i\theta}}{2}\left(\frac{\partial}{\partial r}-i\frac{1}{r}\frac{\partial}{\partial\theta}\right)
\\
\\
\displaystyle
\frac{\partial}{\partial x_2}
=\frac{-ie^{i\theta}}{2}\left(\frac{\partial}{\partial r}+i\frac{1}{r}\frac{\partial}{\partial\theta}\right)
+\frac{ie^{-i\theta}}{2}\left(\frac{\partial}{\partial r}-i\frac{1}{r}\frac{\partial}{\partial\theta}\right)
\end{array}
$$
give the estimate
$$\begin{array}{c}
\displaystyle
\vert\frac{\partial}{\partial x_1}S_N(x;s)\vert+\vert\frac{\partial}{\partial x_2}S_N(x;s)\vert
\le
s
\frac{d}{dz}\{E_{1/n}(z)-E^{nN}_{1/n}(z)\}\vert_{z=Rs}\\
\\
\displaystyle
+k^2R\{E_{1/n}(z)-E^{nN}_{1/n}(z)\}\vert_{z=Rs}.
\end{array}
\tag {3.10}
$$
From the proof of Proposition 3.2 in \cite{IS2} one has
$$\displaystyle
\vert E_{1/n}(z)-E_{1/n}^{nN}(z)\vert
\le\sum_{l=1}^n\frac{\vert z\vert^{nN+l}}{\displaystyle\Gamma\left(N+1+\frac{l}{n}\right)}
\,e^{\displaystyle\vert\text{Re}\,z^n\vert}
\tag {3.11}
$$
and
$$\displaystyle
\vert\frac{d}{dz}\{E_{1/n}(z)-E_{1/n}^{nN}(z)\}\vert
\le n\vert z\vert^{n-1}
\sum_{l=1}^n\frac{\vert z\vert^{n(N-1)+l}}
{\displaystyle\Gamma\left(N+\frac{l}{n}\right)}\,e^{\displaystyle\vert\text{Re}\,z^n\vert}.
\tag {3.12}
$$
Note that these are sharper than (3.5) and (3.6) of Proposition 3.2 in \cite{IS2}.
Consider the case when $z=Rs(N)$.  Then we get, as $N\longrightarrow\infty$
$$\begin{array}{c}
\displaystyle
\sum_{l=1}^n\frac{\vert z\vert^{nN+l}}{\displaystyle\Gamma\left(N+1+\frac{l}{n}\right)}
\,e^{\displaystyle\vert\text{Re}\,z^n\vert}
=O\left(\frac{(Rs(N))^{n(N+1)}e^{(Rs(N))^n}}{\displaystyle\Gamma\left(N+1+\frac{1}{n}\right)}\right);
\\
\\
\displaystyle
s(N)n\vert z\vert^{n-1}
\sum_{l=1}^n\frac{\vert z\vert^{n(N-1)+l}}
{\displaystyle\Gamma\left(N+\frac{l}{n}\right)}\,e^{\displaystyle\vert\text{Re}\,z^n\vert}
=O\left(\frac{(Rs(N))^{n(N+1)}e^{(Rs(N))^n}}{\displaystyle
\Gamma\left(N+\frac{1}{n}\right)}\right).
\end{array}
\tag {3.13}
$$
Since
$$\displaystyle
\frac{1}{\displaystyle\Gamma\left(N+\frac{1}{n}\right)}=O\left(\frac{1}{(N-1)!}\right),
$$
from (3.9), (3.10) and (3.11) to (3.13) we obtain
$$\begin{array}{c}
\displaystyle\vert S_N(x;s(N))\vert+
\vert\frac{\partial}{\partial x_1}S_N(x;s(N))\vert+\vert\frac{\partial}{\partial x_2}S_N(x;s(N))\vert
\\
\\
\displaystyle
=O\left(\frac{N^2(Rs(N))^{n(N-1)}e^{(Rs(N))^n}}{(N-1)!}\right).
\end{array}
\tag {3.14}
$$
Using the Stirling formula, we have
$$\displaystyle
\frac{\xi(N)^{N-1}e^{\xi(N)}}
{(N-1)!}
=O\left(N^{-1/2}e^{\displaystyle N\left(\frac{\gamma}{e}+\log\gamma\right)}\right)
\tag {3.15}
$$
where $\{\xi(N)\}_{N=1,\cdots}$ is an arbitrary sequence of positive numbers satisfying, as $N\longrightarrow\infty$
$\displaystyle
\xi(N)=(\gamma/e)N+O(1)$
and $0<\gamma<\gamma_0$.  Now the conclusion follows from (3.7), (3.14) and (3.15).

\noindent
$\Box$

From (3.5) we know that $g^{1/n}_{(y,\,\omega)}(\,\cdot\,;s,k)_N$ given by (1.1) has the expression
$$\displaystyle
g^{1/n}_{(y,\,\omega)}(\varphi;s,k)_N
=e^{-iky\cdot\varphi}g_N(\varphi;s,k,\omega),\,\,\varphi\in S^1.
$$
Then from Definition 3.1, (3.8) and the equation $Hg^{1/n}_{(y,\omega)}(\,\cdot\,;s,k)_N(x)=Hg_N(\,\cdot\,;s,k,\omega)(x-y)$,
we immediately obtain

\proclaim{\noindent Corollary 3.1.}  Let $\{s(N)\}_{N=1,\cdots}$
be same as in Theorem 3.1.
Then for any fixed $(y,\omega)\in B_R\times S^1$ we have, as $N\longrightarrow\infty$
$$\begin{array}{c}
\displaystyle \sup_{\vert x\vert\le
R}\vert Hg_{(y,\,\omega)}^{1/n}(\,\cdot\,;s(N),k)_N(x)
-E_{1/n}(x-y;s(N), k,\omega)\vert
\\
\\
\displaystyle
+\sup_{\vert x\vert\le R}\vert\nabla
\{Hg_{(y,\,\omega)}^{1/n}(\,\cdot\,;s(N),k)_N(x)
-E_{1/n}(x-y;s(N),k,\omega)\}\vert
=O(N^{-\infty}).
\end{array}
\tag {3.16}
$$
\endproclaim

\section{Remarks}

{\bf\noindent Remark 4.1.}
It should be pointed out that the density (1.1) satisfies
$$\displaystyle
\left(\Phi_x,\,g^{1/n}_{(x,\,\omega)}(\,\cdot\,;s,k)_N\right)_{L^2(S^1)}=1.
$$
Therefore we obtain the relationship between our indicator function and Kirsch's one:
$$\displaystyle
K(x)\le\vert I^{1/n}_{(x,\,\omega)}(s(N))_N\vert.
$$
This together with Theorem 1.1 explains why $K(x)=0$ in a case
that the point $x$ can be connected with infinity by a straight
line without intersecting $\overline D$.  Note that, if $x\in\,D$,
then $0<K(x)<\infty$, however, from Theorem 1.1 we know that, for
all $\omega\in\,S^1$ $\lim_{N\longrightarrow\infty}\vert
I^{1/n}_{(x,\,\omega)}(s(N))_N\vert=\infty$.

{\bf\noindent Remark 4.2.}
Theorem 1.1 does not cover the `critical' case when
both $\overline{C_{y}(\omega,\pi/2n)}\cap\overline D\not=\emptyset$
and $C_{y}(\omega,\pi/2n)\cap D=\emptyset$ are satisfied.
This is coming from a lack of a necessary uniform estimate
of the function $E_{\alpha}(x;s,k,\omega)$.
At the present time we do not know what one can say about the behaviour of the indicator
function as $N\longrightarrow\infty$ in this case.
Note that the results in \cite{Ik9} and Theorem 1.1
in \cite{IS2} completely cover this type case.

\centerline{{\bf Acknowledgement}}

This research was partially supported by Grant-in-Aid for
Scientific Research (C)(No. 21540162) of Japan Society for
the Promotion of Science.

\section{Appendix.  Proof of Lemma 2.1}

\noindent
Write
$$
\displaystyle
u(x;\tau)-f(\tau(x_1+ix_2))
=-\left(\frac{k\vert x\vert}{2}\right)^2\hat{J}_1(k\vert x\vert)
\int_0^1 f(\tau t(x_1+ix_2))dt
+\left(\frac{k\vert x\vert}{2}\right)^2R(x;\tau)
$$
where
$$\displaystyle
R(x;\tau)=\int_0^1f(\tau t(x_1+ix_2))\{\hat{J}_1(k\vert x\vert)
-\hat{J}_1(k\vert x\vert\sqrt{1-t})\}dt.
$$
Using the expression
$$\displaystyle
\hat{J}_m(k\vert x\vert\sqrt{1-t})=m!\sum_{n=0}^{\infty}\frac{\displaystyle (-1)^n}
{\displaystyle (m+n)!n!}\left(\frac{k\vert x\vert}{2}\right)^n(1-t)^n,
\tag {A.1}
$$
one obtains
$$\displaystyle
\frac{d}{dt}\{\hat{J}_1(k\vert x\vert\sqrt{1-t})\}
=\frac{1}{2}\left(\frac{k\vert x\vert}{2}\right)^2\hat{J}_2(k\vert x\vert\sqrt{1-t}).
$$
Then the mean value theorem yields
$$
\vert\hat{J}_1(k\vert x\vert\sqrt{1-t})-\hat{J}_1(k\vert x\vert)\vert
\le\frac{1}{2}\left(\frac{k\vert x\vert}{2}\right)^2t.
$$
From this one gets
$$
\vert R(x;\tau)\vert\le\frac{1}{2}
\left(\frac{k\vert x\vert}{2}\right)^2\int_0^1\vert f(\tau t(x_1+ix_2))\vert tdt.
$$
Now (2.4) and (2.5) are clear.

\noindent
Next from (A.1) we have
$$
\frac{\partial}{\partial x_j}\{\hat{J}_1(k\vert x\vert\sqrt{1-t})\}
=-\left(\frac{k}{2}\right)^2x_j(1-t)\hat{J}_2(k\vert x\vert\sqrt{1-t})
$$
and this yields
$$
\begin{array}{lr}
\displaystyle
\frac{\partial}{\partial x_j}\{u(x;\tau)-f(\tau(x_1+ix_2))\}
=-\left(\frac{k}{2}\right)^22x_j\int_0^1f(\tau t(x_1+ix_2))\hat{J}_1(k\vert x\vert\sqrt{1-t})dt\\
\\
\displaystyle
-\left(\frac{k\vert x\vert}{2}\right)^2
\tau i^{j-1}\int_0^1tf'(\tau t(x_1+ix_2))\hat{J}_1(k\vert x\vert\sqrt{1-t})dt\\
\\
\displaystyle
+\left(\frac{k\vert x\vert}{2}\right)^2\left(\frac{k}{2}\right)^2x_j
\int_0^1f(\tau t(x_1+ix_2))(1-t)\hat{J}_2(k\vert x\vert\sqrt{1-t})dt\\
\\
\displaystyle
=\left(\frac{k\vert x\vert}{2}\right)^2\{A_j(x;\tau)+R_j(x;\tau)\}
\end{array}
\tag {A.2}
$$
where
$$\begin{array}{lr}
\displaystyle
A_j(x;\tau)
=-\frac{2x_j}{\displaystyle\vert x\vert^2}\hat{J}_1(k\vert x\vert)
\int_0^1f(\tau t(x_1+ix_2))dt
-\tau i^{j-1}\hat{J}_1(k\vert x\vert)
\int_0^1 tf'(\tau t(x_1+ix_2))dt\\
\\
\displaystyle
+\left(\frac{k}{2}\right)^2x_j\hat{J}_2(k\vert x\vert)
\int_0^1f(\tau t(x_1+ix_2))dt
\end{array}
$$
and
$$\begin{array}{lr}
\displaystyle
R_j(x;\tau)
=-\frac{\displaystyle 2x_j}{\displaystyle \vert x\vert^2}\int_0^1f(\tau t(x_1+ix_2))
\{\hat{J}_1(k\vert x\vert\sqrt{1-t})-\hat{J}_1(k\vert x\vert)\}dt\\
\\
\displaystyle
-\tau i^{j-1}\int_0^1 tf'(\tau t(x_1+ix_2))
\{\hat{J}_1(k\vert x\vert\sqrt{1-t})-\hat{J}_1(k\vert x\vert)\}dt\\
\\
\displaystyle
+\left(\frac{k}{2}\right)^2x_j
\int_0^1f(\tau t(x_1+ix_2))
\{(1-t)\hat{J}_2(k\vert x\vert\sqrt{1-t})-\hat{J}_2(k\vert x\vert)\}dt.
\end{array}
\tag {A.3}
$$
Change of variables and integration by parts yield
$$\begin{array}{lr}\displaystyle
A_j(x;\tau)
=-\frac{\displaystyle i^{j-1}\hat{J}_1(k\vert x\vert)}{x_1+ix_2}f(\tau(x_1+ix_2))\\
\\
\displaystyle
+\left\{-\frac{\displaystyle (-i)^{j-1}\hat{J}_1(k\vert x\vert)}
{\displaystyle x_1-ix_2}
+\left(\frac{k}{2}\right)^2x_j\hat{J}_2(k\vert x\vert)\right\}\frac{1}{\tau}
\int_0^{\tau}f(w(x_1+ix_2))dw.
\end{array}
\tag {A.4}
$$
Since
$$\displaystyle
\frac{d}{dt}\{\hat{J}_2(k\vert x\vert\sqrt{1-t})\}
=\frac{1}{3}\left(\frac{k\vert x\vert}{2}\right)^2\hat{J}_3(k\vert x\vert\sqrt{1-t}),
$$
one knows that
$$
\vert (1-t)\hat{J}_2(k\vert x\vert\sqrt{1-t})-\hat{J}_2(k\vert x\vert)\vert
\le \left\{\frac{1}{3}\left(\frac{k\vert x\vert}{2}\right)^2+1\right\}t.
$$
Using this together with (A.2), (A.3) and (A.4), we obtain (2.6) and (2.7).

$\Box$

\vskip1cm
\noindent
e-mail address

ikehata@math.sci.gunma-u.ac.jp
\end{document}